\documentclass[12pt]{article}
\usepackage{amsmath, amsthm, amsfonts, amssymb}

\theoremstyle{plain}
\newtheorem{theorem}{Theorem}[section]

\newtheorem{lemma}[theorem]{Lemma}

\theoremstyle{definition}

\newtheorem{definition}[theorem]{Definition}

\numberwithin{equation}{section}

\def \R{\mbox{${\mathbb R}$}}

\begin{document}

\title{Metrizable linear connections in a Lie algebroid}
\author{{\bf Mihai ANASTASIEI} }

\date{Dedicated to Prof. Dr. Constantin Udri\c ste at his 70th anniversary}

\maketitle

\begin{abstract}
A linear connection $D$ in a Lie algebroid is said to be
metrizable if there exists a Riemannian metric $h$ in the Lie
algebroid such that $Dh=0$. Conditions for the linear connection
$D$ to be metrizable are investigated.
\end{abstract}



\section{Introduction}
Le algebroids as particular anchored vector bundles \cite{20} have
now an important place in differential geometry and algebraic
geometry. Initially defined as infinitesimal part of Lie grupoids,
their algebra \cite{12} and geometry is independently and largely
developed \cite{5,8,10,21}. Besides they have proved to be useful
in Mechanics \cite{24,16,2,4,7}, in the theory of nonholonomic
systems\cite{18,9,3}, in control theory \cite{6}, in field theory
\cite{15}, in quantum and classical gravity \cite{22,23}. The
cohomology of Lie algebroids started in \cite{17}. There the
adapted exterior differential $d$ was for the first time
introduced. See also \cite{13}. Holonomy and characteristic
classes have been studied in \cite{11}. For more references see
the monograph \cite{14}

Let $A$ be a Lie algebroid and $D$ an $A-$connection (defined in
the Section 3 ) in a vector bundle $\left(F,q,M\right)$. We say
that $D$ is metrizable if there exists a Riemannian metric $h$ in
$\left(F,q,M\right)$ such that $Dh=0$. The tangent bundle $TM$ is
a trivial Lie algebroid and a $TM-$connection in
$\left(F,q,M\right)$ is nothing but an usual linear connection in
this vector bundle. In \cite{1} we provided  conditions for
metrizability of a $TM-$connection in any vector bundle as well as
in vector bundles endowed with Finsler functions. In this paper we
extend some results from \cite{1} to $A-$connection in
$\left(F,q,M\right)$.The notations from \cite{19}are used.

\section{Lie algebroids}

Let $\xi=\left(E,q,M\right)$ be a vector bundle of rank $m$. Here
$E$ and $M$ are smooth i.e. $C^{\infty}$ manifolds with dim$M=n$,
dim$E=n+m$ and $p:E\longrightarrow M$ is a smooth submersion. The
fibres $E_x=p^{-1}(x)$, $x\in M$ are linear spaces of dimension
$m$ which are isomorphic with the type fibre $\R^m$.

Let $\mathcal F(M)$ be the ring of smooth real functions on $M$.
We denote by $\Gamma(E)$ and $\mathcal X(M)$ the $\mathcal
F(M)-$module of sections of $\xi$ and of the tangent bundle
$\left(TM,\tau,M\right)$, respectively. If
$\left(U,\left(x^i\right)\right)$, $i,j,k\ldots=1,2,\ldots,n$ is a
local chart on $M$, then $\left(\dfrac{\partial}{\partial
x^i}\right)$ provide a local basis for $\mathcal X(U)$. Let
$s_a:U\longrightarrow p^{-1}(U)$, $a,b,c,\ldots=1,2,\ldots,m$ be a
local basis for $\Gamma\left(p^{-1}(U)\right)$. Any section $s$
over $U$ has the form $s=y^as_a(x)$, $x\in U$ and we will take
$\left(x^i,y^a\right)$ as local coordinates on $p^{-1}(U)$. A
change of these coordinates $\left(x^i,y^a\right)\longrightarrow
\left(\tilde{x}^i,\tilde{y}^a\right)$ has the form
\begin{equation}
\begin{split}
\tilde{x}^i&=\tilde{x}^i\left(x^1,\ldots,x^n\right),\quad rank\left(\dfrac{\partial \tilde{x}^i}{\partial x^j}\right)=n,
\\
\tilde{y}^a&=M^a_b(x)y^b,\quad\qquad rank\left(M^a_b(x)\right)=m.
\end{split}
\end{equation}

Let $\xi ^*=\left(E^*,p^*,M\right)$ be the dual of vector bundle
$\xi $ and $\theta^a:U\longrightarrow p^{\ast-1}(U)$,
$x\longrightarrow \theta^a(x)\in E^*_x$ a local basis for
$\Gamma\left(p^{\ast-1}(U)\right)$ such that
$\theta^a(s_b)=\delta^a_b$.

Next, we may consider the tensor bundle of type $(r,s)$, $\mathcal
T^r_s(E)$ over $M$ and its
sections. For $g\in \Gamma\left(E^*\otimes E^*\right)$ we have
$g=g_{ab}(x)\theta^a\otimes \theta^b$. As $\left(E^*\otimes
E^*\right)\cong L_2\left(E,\R\right)$, we may regard $g$ as a
smooth mapping $x\longrightarrow g(x):E_x\times E_x\longrightarrow
\R$ with $g(x)$ a bilinear mapping given by
$g(x)(s_a,s_b)=g_{ab}(x)$, $x\in M$.

If the mapping $g(x)$ is symmetric i.e. $g_{ab}=g_{ba}$ and
positive definite i.e. $g_{ab}(x)\xi ^a\xi ^b>0$ for every $(\xi
^a)\neq0$, one says that $g$ defines a Riemannian metric in the
vector bundle $\xi$.

Let us assume that
\begin{enumerate}
\item[(i)] $\Gamma(E)$ is endowed with a Lie algebra structure
$[,]$ over $\R$, \item[(ii)] There exists a bundle map
$\rho:E\longrightarrow TM$, called anchor map. It  induces a Lie
 algebra homomorphism (denoted also by $\rho$) from $\Gamma(E)$ to $\mathcal X(M)$,
\item[(iii)] For any sections $s_1,s_2\in \Gamma(E)$ and for any $f\in \mathcal F(M)$ the following identity holds
$$
[s_1,fs_2]=f[s_1,s_2]+\rho (s_1)fs_2.
$$
\end{enumerate}

\begin{definition}
The triplet $A=(\xi ,[,],\rho)$ with the properties (i), (ii) and (iii) is called a Lie algebroid.
\end{definition}

\textbf{Examples:}
\begin{enumerate}
\item[1.] The tangent bundle $(TM,\tau ,M)$with the usual Lie bracket and $\rho$ equal to the identity map form
 a Lie algebroid.
\item[2.] Any integrable subbundle of $TM$ with the Lie bracket defined by restriction and $\rho$ the inclusion map
 is a Lie algebroid.
\item[3.] Let $(F,q,M)$ be any vector bundle. On $F$ we have the vertical distribution
 $u\longrightarrow V_uF=$Ker$q_{\ast,u}$, $u\in F$, where $q_{\ast}$ denotes the differential of $q$.

This distribution is integrable. If we regard it as a subbundle of
$TF$, accordingly to Example 2 a Lie algebroid is obtained.

\end{enumerate}

Locally, we set
\begin{equation}
\rho(s_a)=\rho_a^i(x)\dfrac{\partial}{\partial x^i}, \quad [s_a,s_b]=L^c_{ab}(x)s_c.
\end{equation}

The structure function $\rho_a^i$ and $L^c_{ab}$ of the Lie algebroid $A$ have to satisfy the following identities
\begin{equation}
\begin{split}
\rho_a^i\dfrac{\partial \rho_b^j}{\partial
x^i}-\rho_b^i\dfrac{\partial \rho_a^j}{\partial
x^i}=\rho_c^jL^c_{ab},
 \qquad L^c_{ab}+L^c_{ba}=0,\\
\sum_{cycl(abc)}\left(L^d_{ab}L^e_{dc}+\rho_c^i\dfrac{\partial L_{ab}^e}{\partial x^i}\right)=0.
\end{split}
\end{equation}

\section{Connections in Lie algebroids}

Let $A=(\xi ,[,],\rho )$ be a Lie algebroid with $\xi=(E,p,M)$ and let $(F,q,M)$ be any vector bundle.

\begin{definition}
An $A-$connection in the bundle $(F,q,M)$ is a mapping $D:\Gamma(E)\times\Gamma(F)\longrightarrow \Gamma(F)$, $(s,\sigma)\longrightarrow D_s\sigma$ with the properties:
\begin{enumerate}
\item[1)] $D_{s_1+s_2}\sigma=D_{s_1}\sigma+D_{s_2}\sigma$,
\item[2)] $D_{fs}\sigma=fD_s\sigma$,
\item[3)] $D_s{(\sigma _1+\sigma _2)}=D_s{\sigma _1}+D_s{\sigma _2}$,
\item[4)] $D_s{(f\sigma)}=\rho(s)f\sigma+fD_s\sigma $,

for $s,s_1,s_2\in \Gamma(E)$, $\sigma, \sigma_1, \sigma_2\in \Gamma(F)$, $f\in \mathcal F(M)$.
\end{enumerate}
\end{definition}

Notice that a $TM-$connetion in the vector bundle $(F,q,M)$ is nothing but a linear connection in this vector bundle.

\begin{definition}
An $A-$connection in the bundle $\xi =(E,p,M)$ is called a linear connection in the Lie algebroid $A$.
\end{definition}

The notion of tangent lift of a curve on $M$ is generalized as follows.

\begin{definition}
Let $A=(\xi ,[,],\rho)$ be a Lie algebroid with $\xi =(E,p,M)$. A
curve $\alpha :[0,1]\longrightarrow E$ is called admissible or an
$A-$path if $\rho(\alpha(t))=\dfrac{d}{dt}p(\alpha (t))$, $t\in
[0,1]$. The curve $\gamma (t)=p(\alpha(t))$ will be called the
base path of $\alpha$. The $A-$path $\alpha$ is called vertical if
$\rho(\alpha(t))=0$. In this case $\gamma$ reduces to a point and
the curve $\alpha$ is contained in the fibre in that point.
\end{definition}

Locally, if $\alpha(t)=\left(x^i(t),y^a(t)\right)$, then $\gamma (t)=(x^i(t))$ and $\alpha $ is an $A-$path if and only if
\begin{equation}
\rho_a^i(x(t))y^a(t)=\dfrac{dx^i(t)}{dt}, \qquad t\in [0,1],
\end{equation}
and it is a vertical $A-$path if and only if
\begin{equation}
\rho_a^i(x(t))y^a(t)=0, \qquad t\in [0,1].
\end{equation}

Let $(\sigma _{\alpha})$, $\alpha,\beta,\gamma,\ldots=k:=$rank of $(F,q,M)$, a local basis in $\Gamma (F)$.
Then a local section $\sigma$ has the form $\sigma=z^{\alpha}\sigma_{\alpha}$ and $(z^{\alpha})$ are the coordinates in
the fibres of $(F,q,M)$.

For $s=y^as_a$ and $\sigma=z^{\alpha}\sigma_{\alpha}$, by the
Definition 3.1 we have $D_s\sigma=y^a\left(\rho^i_a\dfrac{\partial
z^{\alpha}}{\partial
x^i}+z^{\alpha}D_{s_a})\sigma_{\alpha}\right)$ and if we put
\begin{equation}
D_{s_a}\sigma_{\alpha}=\Gamma^{\beta}_{\alpha a}\sigma_{\beta},
\end{equation}
we get
\begin{equation}\label{1}
D_s\sigma=y^a(D_az^{\beta})\sigma_{\beta},\qquad D_az^{\beta}=\rho^i_a\dfrac{\partial z^{\beta}}{\partial x^i}+\Gamma^{\beta}_{\alpha a}z^{\alpha}.
\end{equation}
For a linear connection $\mathcal D$ in the Lie algebroid $A=(\xi ,[,],\rho)$ we get
\begin{equation}
\mathcal D_s\sigma=y^a(\mathcal D_az^b)s_b,\qquad \mathcal D_az^b=\rho^i_a\dfrac{\partial z^b}{\partial x^i}+\Gamma^b_{ca}z^c.
\end{equation}

Let $D$ be an $A-$connection in the vector bundle $(F,q,M)$ and
$\alpha : [0,1]\longrightarrow E$ an $A-$path.

A smooth mapping $\sigma :[0,1]\longrightarrow F$ is called an $\alpha -$section if $q(\sigma (t))=p(\alpha (t))$, $t\in [0,1]$.

Locally, if $\alpha (t)=(x^i(t),y^a(t))$ then $\sigma (t)=(x^i(t),z^\alpha(t))$.

Let $\Gamma(F)^\alpha $ be the linear space of $\alpha -$section in the vector bundle $(F,q,M)$.
We define an operator $D^{\alpha}:\Gamma(F)^\alpha \longrightarrow \Gamma(F)^\alpha$, $\sigma (t)\longrightarrow
(D^{\alpha}\sigma)(t)$ by
\begin{equation}\label{2}
(D^{\alpha}\sigma)(t)=\left(\dfrac{dz^\beta}{dt}+\Gamma^{\beta}_{\alpha
a}(x(t))z^{\alpha}(t)y^a(t)\right)\sigma_{\beta},
\end{equation}
whenever $\sigma(t)=z^{\alpha}(t)\sigma_{\alpha}$.

The operator $D^{\alpha}$ has the following properties:
\begin{enumerate}
\item[(i)] $D^{\alpha}(c_1\sigma _1+c_2\sigma _2)=c_1D^{\alpha}\sigma_1+c_2D^{\alpha}\sigma_2$,\quad $c_1,c_2\in \R$, $\sigma _1,\sigma _2\in \Gamma(F)^\alpha$,
\item[(ii)] $D^\alpha(fs)=\dfrac{df}{dt}\sigma+fD^\alpha \sigma$, for $\sigma\in \Gamma(F)^\alpha$ and $f:[0,1]\longrightarrow \R$ a smooth function,
\item[(iii)] If $\widetilde{\sigma}$ is a local section that extends $\sigma\in \Gamma(F)^\alpha$ and $\rho(\alpha(t))\neq 0$ ($\rho(\alpha (t))=0$) then $(D^{\alpha}\sigma)(t)=D_{\alpha(t)}\widetilde{\sigma}$ (resp. $(D^{\alpha}\sigma)(t)=D_{\alpha(t)}\widetilde{\sigma}+\dfrac{d\sigma}{dt}$).
\end{enumerate}

The first two properties are immediate by \eqref{2}. To prove
(iii) one uses \eqref{2}, (3.2)  and \eqref{1}

By contradiction one proves that $D^\alpha$ is the unique operator
with the properties (i), (ii) and (iii). Indeed, if
$\widetilde{D}^\alpha :\Gamma(F)^\alpha\longrightarrow
\Gamma(F)^\alpha$ is another operator satisfying (i), (ii) and
(iii), it easily follows that it is has the form that appears in
the second hand of \eqref{2}.

\begin{definition}
An $\alpha-$section $\sigma$ in the vector bundle $(F,q,M)$ is said to be parallel if $D^\alpha\sigma=0$.
\end{definition}

Locally, the $\alpha -$section $\sigma(t)=z^\alpha(t)\sigma_\alpha$ with $\alpha(t)=y^a(t)s_a$ is parallel if and only if
the functions $z^\alpha(t)$ are solutions of the following linear system of ordinary differential equations

\begin{equation}\label{3}
\dfrac{dz^\beta}{dt}+\Gamma^{\beta}_{\alpha a}(x(t))z^{\alpha}(t)y^a(t)=0.
\end{equation}
This system has an unique solution $t\longrightarrow \sigma(t)$
with the initial condition $\sigma(0)=\sigma_0$. This fact allow
us to define the parallel displacement along $\alpha$, denoted by
$P^t_\alpha :F_{\gamma(0)}\longrightarrow
F_{\gamma(t)},\gamma(t)=p(\alpha(t)=q(\sigma(t)), P^t_\alpha
(\sigma_0)=\sigma(t)$. The maps $P^t_\alpha$ are linear
isomorphisms.

In particular, we may take $\alpha$ a loop based at $x\in M$ i.e.
$\gamma(0) =\gamma(1) = x$ and we get the linear isomorphism
$P_{\alpha} :F_x\mapsto F_x$. Its inverse is $P_{\alpha^{-1}}$
where $\alpha^{-1}$ is the reverse loop of $\alpha$ and if we
consider the composite $\alpha_1\circ\alpha_2 $ that is $\alpha_2$
followed by $\alpha_1$ of two loops based on $x$ it comes out that
$P_{\alpha_1\circ\alpha_2}= P_{\alpha_2}\circ P_{\alpha_1} $. On
this way one obtains a subgroup of the linear isomorphisms of
$F_x$ called the holonomy group of $D$, denoted by $\Phi(x))$.

We fix $t$ and consider
$(P^t_\alpha)^{-1}:F_{\gamma(t)}\longrightarrow F_{\gamma(0)}$.
Locally, if
$(P^{\tau}t_\alpha)^{-1}(\sigma(\tau))=\widetilde{z}^\beta(\tau)\sigma_\beta$,
then $\widetilde{z}^\beta(t)=z^\beta(t)$,
$\widetilde{z}^\beta(0)=(P^t_\alpha)^{-1}(\sigma(t))$ and
$(\widetilde{z}^\beta)$ are solutions of \eqref{3}. By Taylor's
formula
$\widetilde{z}^\beta(t)=\widetilde{z}^\beta(0)+t\dfrac{d\widetilde{z}^\beta}{dt}(0)\ldots$,
hence
$z^\beta(t)-z^\beta(0)=((P^t_\alpha)^{-1}(\sigma(t))^\beta-z^\beta(0)+t\dfrac{d\widetilde{z}^\beta}{dt}(0)\ldots$.
We divide this by $t$, take $t\rightarrow 0$ and obtain
\begin{equation}
(D^\alpha\sigma)(0)=lim_{t\rightarrow0}\dfrac{(P^t_\alpha)^{-1}(\sigma(t))-\sigma(0)}{t}.
\end{equation}
Suppose now that $\cal{D}$ is a linear connection in the Lie
algebroid $A=(\xi,[,],\rho)$, $\xi=(E,p,M)$. An $\alpha-$path is
called geodesic if $\cal{D}^\alpha\alpha=0$. Locally, if
$\alpha(t)=(x^i(t),y^a(t))$, then
$\cal{D}^{\alpha}\alpha= [\dfrac{dy^a}{dt}+\Gamma^a_{bc}(x(t))y^b(t)y^c(t) ]s_a$
and $\alpha$ is a geodesic if and only if the functions
$x^i(t),y^a(t)$ are solution of the following system of ordinary
differential equations:
\begin{equation}
\begin{split}
\dfrac{dx^i}{dt}=\rho^i_a(x(t))y^a(t),\\
\dfrac{dy^a}{dt}+\Gamma^a_{bc}(x(t))y^b(t)y^c(t)=0.
\end{split}
\end{equation}
It is clear that one has the existence and uniqueness of geodesic
with a given base point $x\in M$ and a given $y_0\in E_{x_0}$.If
for a pair $(x_0,y_0)$ we have $\rho^i_a(x_0)y^a_0=0$, the
corresponding geodesic is contained in the fibre $E_{x_0}$ i.e. it
is a vertical $A$-path.

\section{Riemannian metrics in Lie algebroids}

Let $A=(\xi,[,],\rho)$ be a Lie algebroid with $\xi=(E,p,M)$ and a
vector bundle $(F,q,M)$ endowed with an $A-$ connection $D$ whose
local coefficients are $(\Gamma^\alpha_{\beta a})$.

A Riemannian metric in $(F,q,M)$ is a mapping $g$ that assigns to
any $x\in M$ a scalar product $g_x$ in $E_x$ such that for any
local section $\sigma_1,\sigma_2\in \Gamma(F)$, the function
$x\longrightarrow g_x(\sigma_1,\sigma_2)$ is smooth. Locally, we
set $g_x(\sigma_\alpha,\sigma_\beta)=g_{\alpha\beta}(x)$ and so
$g_x(\sigma_1,\sigma_2)=g_{\alpha\beta}(x)z^\alpha_1z^\beta_2$ if
$\sigma_1=z^\alpha_1\sigma_\alpha$,
$\sigma_2=z^\beta_2\sigma_\beta$.

The operator of covariant derivative $D$ can be extended to the
tensor algebra of $(F,q,M)$ taking $D_\sigma f=\rho(\sigma)f$,
assuming that it commutes with the contractions and behaves like a
derivation with respect to tensor product. It comes out that if
$\omega$ is a section in the dual $(F^*,q^*,M)$ then
$$
(D_s\omega)(\sigma)=\rho(s)\omega(\sigma)-\omega(D_s\sigma), \quad
s\in \Gamma(E), \quad \sigma\in \Gamma(F)
$$
and if $a$ is a section in $L^2(F,\R)$, then
\begin{equation}\label{4}
(D_sa)(\sigma_1,\sigma_2)=\rho(s)a(\sigma_1,\sigma_2)-a(D_s\sigma_1,\sigma_2)-a(\sigma_1,D_s\sigma_2),
 s\in \Gamma(E), \sigma_1\sigma_2\in \Gamma(F).
\end{equation}

\begin{definition}
We say that the Riemannian metric $g$ is compatible with the
$A-$connection $D$ if $D_sg=0$ for every $s\in \Gamma(E)$.
\end{definition}

By \eqref{4} the condition of compatibility between $g$ and $D$ is
equivalent to
\begin{equation}\label{5}
\rho(s)g(\sigma_1,\sigma_2)=g(D_s\sigma_1,\sigma_2)+g(\sigma_1,D_s\sigma_2),
\quad s\in \Gamma(E), \quad \sigma_1,\sigma_2\in \Gamma(F).
\end{equation}
Locally, \eqref{5} is written as follows
\begin{equation}
\rho_a^i(x)\dfrac{\partial g_{\alpha\beta}}{\partial
x^i}=\Gamma^\gamma_{\alpha
a}(x)g_{\gamma\beta}(x)+\Gamma^\gamma_{\beta
a}(x)g_{\alpha\gamma}(x).
\end{equation}

The operator $D^{\alpha}$ can be extended to $\alpha-$ section in
the tensor bundles constructed with $(F,q,M)$ and one deduces that

\begin{equation}
(D^{\alpha}g(t)(\sigma_1(t),\sigma_2(t))=
(\dfrac{g_{\mu\nu}}{dt}-g_{\mu\eta}\Gamma^{\eta}_{\mu a}y^a
-g_{\eta\nu}\Gamma^{\eta}_{\mu a}y^a)\sigma_1(t)\sigma_2(t)
\end{equation}

If $(F,q,M)$ coincides with $(E,p,M)$ we have

\begin{theorem}\label{T1}
There exists an unique linear connection $\nabla$ in the Lie
algebroid $A$ such that
\begin{enumerate}
\item[(i)] $\nabla_sg=0$,
\item[(ii)] $\nabla_{s_1}s_2-\nabla_{s_2}s_1=[s_1,s_2]$, $s,s_1,s_2\in \Gamma(E)$.
\end{enumerate}
It is given by the formula
\begin{equation}\label{6}
\begin{split}
2g\left(\nabla_{s_1}s_2,s_3\right)=&\rho(s_1)g(s_2,s_3)+\rho(s_2)g(s_1,s_3)-\rho(s_3)g(s_1,s_2)\\
&+g([s_3,s_1],s_2)+g([s_3,s_2],s_1)+g([s_1,s_2],s_3)
\end{split}
\end{equation}
and its local coefficients are given by
\begin{equation}\label{7}
\Gamma^a_{bc}=\dfrac{1}{2}g^{ad}\left(\rho_b^i\dfrac{\partial
g_{cd}}{\partial x^i}+\rho_c^i\dfrac{\partial g_{bd}}{\partial
x^i}-\rho_d^i\dfrac{\partial g_{bc}}{\partial
x^i}+L^e_{dc}g_{eb}+L^e_{db}g_{ec}-L^e_{bc}g_{ed}\right).
\end{equation}
\end{theorem}

\begin{proof}
In the condition (i) written for $s_1,s_2,s_3\in \Gamma(E)$ we
cyclically permute $s_1,s_2,s_3$ and so we obtain two new
identities. We add these and from the result we subtract the
first. Using (ii) some terms cancel each other and we get
\eqref{6}. Writing \eqref{6} in a local basis of sections we find
\eqref{7}. The uniqueness follows by contradiction.
\end{proof}

If we put
\begin{equation}
T_\nabla(s_1,s_2)=\nabla_{s_1}s_2-\nabla_{s_2}s_1-[s_1,s_2], \quad
s_1,s_2\in \Gamma(E)
\end{equation}
we get a section in the bundle $L(E,E;E)$ that may be called the
torsion of $\nabla$.

The curvature of $\nabla$ is defined by
\begin{equation}
R_\nabla(s_1,s_2)s_3=\nabla_{s_1}\nabla_{s_2}s_3-\nabla_{s_2}\nabla_{s_1}s_3-\nabla_{[s_1,s_2]}s_3,
\quad s_1,s_2,s_3\in \Gamma(E).
\end{equation}

The connection $\nabla$ given by the Theorem \ref{T1} is called
the Levi-Civita connection of $A$.

We stress that the Theorem \ref{T1} says that given $g$ there
exists and is unique $\nabla$ such that $\nabla g=0$ and $T_\nabla
=0$. Now we give a different proof of this theorem.

Given $g$ we may associate to it the energy function $\mathcal{E}
:E\longrightarrow \R$, $\mathcal{E} (s)=g(s,s)$, $s\in E$.
Locally, $\mathcal{E}(x,y)=g_{ab}(x)y^ay^b$, $s=y^as_a$.

The energy function $\mathcal{E}$ is a regular Lagrangian on $E$
i.e. $\det\left(\dfrac{1}{2}\dfrac{\partial
^2\mathcal{E}}{\partial y^a\partial
y^b}\right)=\det(g_{ab}(x))\neq0$.

In \cite{2}, we associated to any regular Lagrangian $L$ on a Lie
algebroid a semispray on $E$ i.e. a vector field
$$
S=\rho_a^iy^a\dfrac{\partial}{\partial
x^i}-2G^a_L(x,y)\dfrac{\partial}{\partial y^a}
$$
with
\begin{equation}\label{8}
G^a_L=\dfrac{1}{4}g^{ab}\left(\dfrac{\partial ^2L}{\partial y^b
\partial x^i}\rho_c^iy^c-\rho_b^i\dfrac{\partial L}{\partial
x^i}-L^c_{bd}y^d\dfrac{\partial L}{\partial y^c}\right),
\end{equation}
where $g_{ab}=\dfrac{1}{2}\dfrac{\partial ^2L}{\partial y^a
\partial y^b}$ and $(g^{ab})$ is the inverse of the matrix
$(g_{ab})$.

Taking $L=\mathcal{E}$ in \eqref{8}, a direct calculation in which
$L^a_{cd}y^cy^d=0$ is used, shows that the semispray associated to
$\mathcal{E}$ has the form
\begin{equation}\label{9}
S=\rho_a^iy^a\dfrac{\partial}{\partial
x^i}-\Gamma^a_{cd}(x)y^cy^d\dfrac{\partial}{\partial y^a},
\end{equation}
with $\Gamma^a_{cd}$ given by \eqref{7}. These coefficients
determines $\nabla$. They are symmetric in bottom indices, hence
$T_\nabla =0$. The uniqueness of $\nabla$ follows by
contradiction.

Note that \eqref{9} gives a $2-$homogeneous semispray, that is a
spray.

By \eqref{9} it follows
\begin{theorem}
The integral curves of $S$ are just the geodesics of the
Levi-Civita connection $\nabla$ in the Lie algebroid $A$.
\end{theorem}
For a different derivation of $S$ from $g$ we refer to \cite{24}.

Now we come back to the general framework and prove
\begin{lemma}\label{L1}
Let be the vector bundle $(F,q,M)$ endowed with an
$A-$connec\-ti\-on $D$ and a Riemannian metric $g$. Then for any
$A-$path $\alpha:t\longrightarrow \alpha (t)$, $t\in [0,1]$ on $E$
with base curve $\gamma =q\circ \alpha= p\circ \alpha$ we have
\begin{equation}
\begin{split}
\left(D_{\alpha}g\right)\left(\sigma_1,\sigma_2\right)=\lim_{t\rightarrow
0}\dfrac{1}{t}\left[g_{\gamma(t)}\left(P_\alpha^t\sigma_1,P_\alpha^t\sigma_2\right)-g_{\gamma(0)}\left(\sigma_1,\sigma_2\right)\right],
\end{split}
\end{equation}
where $\sigma_1,\sigma_2\in F_{\gamma(0)}$ and
$P_\alpha^t:F_{\gamma(t)}\longrightarrow F_{\gamma(0)}$ is the
parallel displacement defined by $D$ along $\alpha$.
\end{lemma}
\begin{proof}
Let $\widetilde{\sigma_1}$ and $\widetilde{\sigma_2}$ be the
parallel $\alpha-$sections in $(F,q,M)$ such that
$\widetilde{\sigma_1}(0)=\sigma_1$,
$\widetilde{\sigma_2}(0)=\sigma_2$. Then
$P_\alpha^t\sigma_1=\widetilde{\sigma_1}(t)$ and
$P_\alpha^t\sigma_2=\widetilde{\sigma_2}(t)$. By Taylor formula in
the local basis $(\sigma_\beta)$ we get
$\left(P_\alpha^t\sigma_1\right)^\mu=\sigma_1^\mu-\Gamma^\mu_{\nu
a}(\tau)\widetilde{\sigma_1}^\nu(\tau)y^a(\tau)t$ for $\tau\in
(0,t)$ and a similar formula for $P_\alpha^t\sigma_2$. Recall that
$\alpha (t)=y^a(t)s_a$. Then by using again the Taylor formula and
omitting the terms which contain $t^2$ and $t^3$, we may write
\begin{equation*}
\begin{split}
&g_{\mu\nu}(\gamma(t))\left(P_\alpha^t\sigma_1\right)^\mu \left(P_\alpha^t\sigma_2\right)^\nu-g_{\mu\nu}(\gamma(0))\sigma_1^\mu \sigma_2^\nu=\\
=&\left(g_{\mu\nu}(\gamma(0))+\dfrac{dg_{\mu\nu}}{dt}(\gamma(0))t\right)\left(P_\alpha^t\sigma_1\right)^\mu \left(P_\alpha^t\sigma_2\right)^\nu-g_{\mu\nu}(\gamma(0))\sigma_1^\mu \sigma_2^\nu\\
=&\left(\dfrac{dg_{\mu\nu}}{dt}-g_{\mu\eta}\Gamma^\eta_{\nu
a}-g_{\eta\nu}\Gamma^\eta_{\mu a}\right)\sigma_1^\mu \sigma_2^\nu
t,
\end{split}
\end{equation*}
where the terms in the parenthesis are computed for some values
$\tau$ in $(0,t)$.

Now dividing by $t$, taking $t\rightarrow0$ and looking at (4.4), we obtain (4.10).
\end{proof}

\begin{definition}
An $A-$conection $D$ in $(F,q,M)$ is said to be metrizable if
there exists a Riemannian metric in $(F,q,M)$ such that $Dg=0$.
\end{definition}

Based on  Lemma \ref{L1} we obtain:

\begin{theorem}
Any $A-$conection $D$ is metrizable with respect to $g$ if and
only if all its parallel displacements are isometries with respect
to $g$.
\end{theorem}

In particular, the holonomy group $\Phi(x)$ is made up of
isometries of $(F_x,$ $g(x))$.

Using a local chart around $x$ we may put $\Phi(x)$ in an $1:1$
correspondence with a subgroup $\widetilde{\Phi}(x)$ of
$GL(k,\R)$, $k=$rank$F$. A change of local chart moves that
subgroup in a conjugate of it. We identify $\Phi(x)$ with this
class of conjugate subgroups in $GL(k,\R)$. With this
identification, by Theorem 4.6 we get

\begin{theorem}
A necessary condition for an $A-$connection to be metrizable is
that the holonomy group $\Phi(x)$ be a subgroup of the orthogonal
group $O(k)$ for every $x\in M$.
\end{theorem}
Indeed, if $\Phi(x)$ is made up of isometries of $(F_x,g(x))$, the
elements of $\widetilde{\Phi}(x)$ are isometries of $(\R^k,<,>)$,
with the inner product $<,>$ induced by $g(x)$.

Authors' address:

"Al.I.Cuza" University, Iasi

Faculty of Mathematics,

700506,  Ia\c si, Romania

Email : anastas@uaic.ro


\begin{thebibliography}{99}

\bibitem{1}
Anastasiei, M., {\it Metrizable linear connections in vector
bundles}, Publ.Math.Debrecen,vol.62(2003),277-287.

\bibitem{2}
Anastasiei, M.,{\it Geometry of Lagrangians and semisprays on Lie
algebroids.} Proceedings of the 5th Conference of Balkan Society
of Geometers, 10--17, BSG Proc., 13, Geom. Balkan Press,
Bucharest, 2006.

\bibitem{3}
Anastasiei, M., {\it  Mechanical systems on Lie algebroids.}
Algebras Groups Geom. 23 (2006), no. 3, 235--245.

\bibitem{4}
Anastasiei, Mihai {\it Semisprays on Lie algebroids.
Applications.} Tensor (N.S.) 69 (2008), 190--198.

\bibitem{5}
Boucetta M., {\it Riemannian Geometry of Lie algebroids}
arXiv:0866.3522v1[math.DG],21.06.2008

\bibitem{6}
Cortés, Jorge; de León, Manuel; Marrero, Juan C.; Martín de Diego,
D.; Martínez, Eduardo, {\it  A survey of Lagrangian mechanics and
control on Lie algebroids and groupoids.} Int. J. Geom. Methods
Mod. Phys. 3 (2006), no. 3, 509--558.

\bibitem{7}
Crasmareanu, Mircea; Hretcanu, Cristina-Elena {\it Last
multipliers on Lie algebroids. }Proc. Indian Acad. Sci. Math. Sci.
119 (2009), no. 3, 287--296.

\bibitem{8}
Degeratu, Marian; Ivan, Mihai, {\it  Linear connections on Lie
algebroids.} Proceedings of the 5th Conference of Balkan Society
of Geometers, 44--53, BSG Proc., 13, Geom. Balkan Press,
Bucharest, 2006.

\bibitem{9}
de Leon, Manuel; Marrero, Juan C.; Martínez, Eduardo, {\it
Lagrangian submanifolds and dynamics on Lie algebroids.} J. Phys.
A 38 (2005), no. 24, R241--R308.

\bibitem{10}
Grabowski, J.; Urbanski, P. {\it Tangent and cotangent lifts and
graded Lie algebras associated with Lie algebroids.} Ann. Global
Anal. Geom. 15 (1997), no. 5, 447--486.

\bibitem{11}
Fernandes, Rui Loja {\it Lie algebroids, holonomy and
characteristic classes.} Adv. Math. 170 (2002), no. 1, 119--179.

\bibitem{12}
Higgins, Philip J.; Mackenzie, Kirill, {\it Algebraic
constructions in the category of Lie algebroids.} J. Algebra 129
(1990), no. 1, 194--230.

\bibitem{13}
Ivan, Gheorghe, {\it  Cohomology on Lie algebroids.} Tensor (N.S.)
67 (2006), no. 3, 283--293.

\bibitem{14}
Mackenzie, Kirill C. H., { \it General theory of Lie groupoids and
Lie algebroids.} London Mathematical Society Lecture Note Series,
213. Cambridge University Press, Cambridge, 2005. xxxviii+501 pp.


\bibitem{15}
 Martinez, Eduardo, {\it Classical field theory on Lie algebroids:
 variational aspects.} J. Phys. A 38 (2005), no. 32, 7145--7160.

\bibitem{16}
Martinez, E.,{\it  Lagrangian mechanics on Lie algebroids.} Acta
Appl. Math. 67 (2001), no. 3, 295--320.

\bibitem{17}
 Maxim-Raileanu, Liliana'\,{ \it Cohomology of Lie algebroids.}
An. Sti. Univ. ``Al. I. Cuza'' Iasi Sect. I a Mat. (N.S.) 22
(1976), no. 2, 197--199.

\bibitem{18}
Mestdag, Tom; Langerock, Bavo, {\it A Lie algebroid framework for
non-holonomic systems.} J. Phys. A 38 (2005), no. 5, 1097--1111.

\bibitem{19}
Miron, R., Anastasiei, M., {\it The Geometry of Lagrange Spaces:
Theory and Applications}, Kluwer Academic Publishers. FTPH 59,
1994.

\bibitem{20}
Popescu, Paul; Popescu, Marcela. {\it  Anchored vector bundles and
Lie algebroids.} Lie algebroids and related topics in differential
geometry (Warsaw, 2000), 51--69, Banach Center Publ., 54, Polish
Acad. Sci., Warsaw, 2001.


\bibitem{21}
Popescu, Liviu, {\it Geometrical structures on Lie algebroids.}
Publ. Math. Debrecen 72 (2008), no. 1-2, 95--109.

\bibitem{22}
Vacaru S., {\it Clifford - Finsler algebroids and nonholonomic
Einstein -Dirac structures.}J. Math. Phys. 47(2006),093504

\bibitem{23}
Vacaru S., {\it Einstein  - Cartan  algebroids and black holes in
solitonic background. } arXiv ; gr-qc/0501057,35 p.

\bibitem{24}
Weinstein A., {\it Lagrangian Mechanics and Grupoids},
Fields Institute Communications, vol. 7, 1996, p. 207-231.


\end{thebibliography}
\end{document}